\documentclass{snapshotmfo}

\categorizationmath{geometry and topology} 
\categorizationconnect{physics} 
\license{CC-BY-NC-SA-4.0} 
\snapshotid{8}{2020}
\junioreditor{-}{junior-editors@mfo.de}
\senioreditor[f]{Sophia Jahns (for Carla Cederbaum)}{senior-editor@mfo.de}
\director[m]{Gerhard Huisken}

\usepackage[utf8]{inputenc}
\usepackage{amsmath,amssymb}

\usepackage[english]{babel}

\author{Steven Rayan \thanks{Steven Rayan is supported through the Natural Sciences and Engineering Research Council of Canada Discovery Grant program and the Canadian Tri-Agency New Frontiers in Research Fund (Exploration Stream).} \and Laura P. Schaposnik \thanks{Laura P. Schaposnik is supported by the National Science Foundation through  the grant DMS-1509693, the CAREER Award DMS-1749013, and by the Alexander von Humboldt Foundation. This material is also based upon work supported by the National Science Foundation under Grant No. DMS-1440140 while the author was in residence at the Mathematical Sciences Research Institute in Berkeley, California, during the Fall 2019 semester.}}

\title{Higgs bundles without geometry}

\authorinfo{\authorname{Steven Rayan} is a professor of pure mathematics at the University of Saskatchewan.}
\authorinfo{\authorname{Laura P. Schaposnik} is a professor of pure mathematics at the University of Illinois at Chicago.}

\begin{document}

\begin{abstract}
Higgs bundles appeared a few decades ago as solutions to certain equations from physics and have attracted much attention in geometry as well as other areas of mathematics and physics.  Here, we take a very informal stroll through some aspects of linear algebra that anticipate the deeper structure in the moduli space of Higgs bundles.
\end{abstract}

\section{Introduction}
Higgs bundles have been making waves in mathematics for over 30 years now.  They are solutions to  certain differential equations that originate in mathematical physics --- specifically, the self-dual Yang-Mills equations reduced by two dimensions \cite{Hi1} --- but have become staples in geometry (algebraic, differential, and symplectic), representation theory, and even number theory.  In a spectacular way, Higgs bundles were used to prove the Fundamental Lemma, a Fields Medal-worthy result \cite{ngo}.  Coming full circle, Higgs bundles have a renewed importance in high-energy physics through applications to string theory and mirror symmetry.


\section{Moduli spaces in modern and ancient life}
We will delay the definition of a Higgs bundle until the next section, because
on its own, a single Higgs bundle is not so important.  What makes them special is the confluence of geometric features of an entire family of Higgs bundles, known as the \emph{moduli space}. Let's take a moment to say a few words about moduli spaces before we go into Higgs bundles themselves.

A moduli space is somewhat akin to a telephone book in which you can look up every Higgs bundle (or some other kind of mathematical object).  If we continue with the telephone analogy, note that --- especially these days ---\linebreak a single person might have multiple phone numbers, say, a home landline, a mobile, and a work number.  To form the moduli space of phone numbers, we consider someone's home number, mobile number, and work number to be three interchangeable avatars for the same thing: the person you are trying to reach.  As such, we make the decision that, for any given person, our phone number moduli space only lists one of these numbers.\footnote{The practicality here is debatable.  On the one hand, the book is data efficient because it only contains one entry per person.  On the other hand, it only gives you one way to reach a person.} This is what we call an \emph{equivalence relation} on a set, in this case the set of all telephone numbers: two numbers are deemed equivalent if they belong to the same person.  Our moduli space of Higgs bundles does the same thing.  An individual Higgs bundle might have many avatars.  We can form the moduli space by picking a preferred avatar for each.

If the telephone example seems too impractical, note that the set of time-telling hours is an ancient moduli space of immense practicality.  Humanity opted out of telling one another the number of hours that have elapsed since the dawn of known history.  Instead, we reset the clock every 12 or 24 hours.  Going with the former (the AM/PM clock), we elect to treat 1 and 13 and 25 and so on as avatars of the same hour.   We keep only 1 o'clock and throw away the rest of its avatars, which is why we never refer to 37 o'clock or 49 o'clock, etc.   Mathematicians would use the notation $\mathbb Z_{12}$ to denote the moduli space of time-telling AM/PM hours.

\section{Vector bundles and matrices with polynomial entries}

So what is a Higgs bundle, then?  To the reader who has not been exposed to serious complex geometry --- Riemann surfaces, vector bundles, Jacobians, holomorphic differentials ---  there is not much that can be said precisely.  Accepting some imprecision, a \textit{vector bundle} can be imagined as a surface together with lines or vector spaces at each point on the surface.  

\begin{figure}[ht]
\centering
\resizebox{0.75\textwidth}{!}{%
  \includegraphics{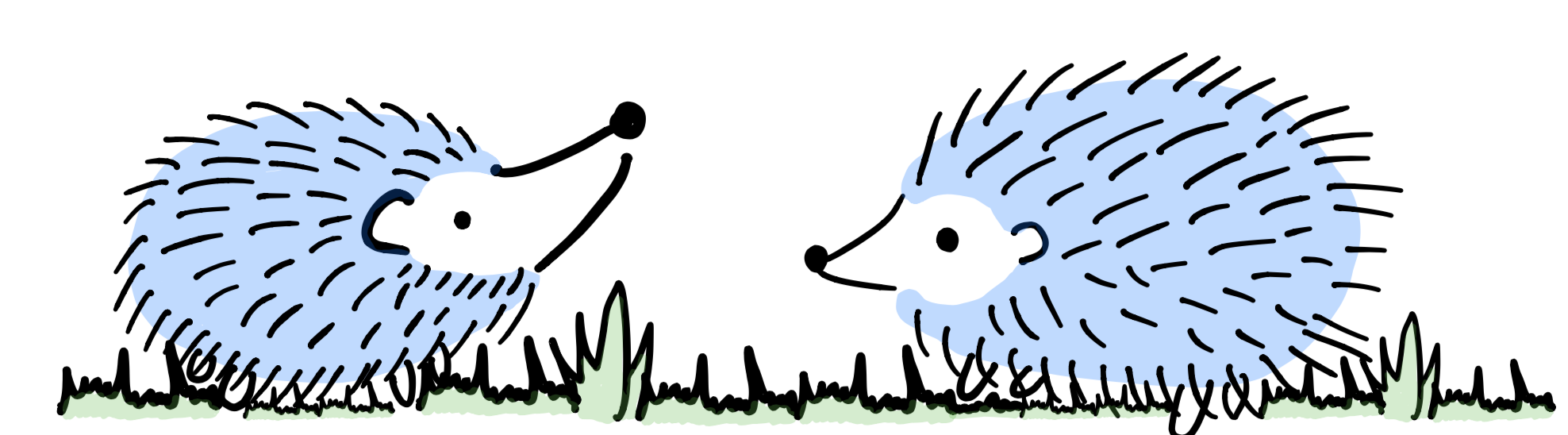}   
  }
\caption{A line bundle over a hedgehog given by its thin hairs. }\label{Lau1}
\end{figure}

In nature, we can think of a hedgehog as a fairly good real-life cartoon for a vector bundle: indeed, over each point of their skin, the hedgehog has a hair which determines a $1$-dimensional space, and thus we would call the hedgehog a (real) \textit{line bundle} as shown in Figure \ref{Lau1}. If we replace each hair with a piece of cardboard, then we will have a rank-2 vector bundle. If we replace each sheet of cardboard with a box, then we will have a rank-3 vector bundle; and so on. We refer to the skin of the hedgehog as the \emph{base space} of these bundles. 

A crucial property of vector bundles is that of \emph{nontriviality}.  This allows a bundle to behave differently over different small patches of the base space, and is useful in many applications --- particularly in physics.  A good base space will have patches that are each in correspondence with $\mathbb R^n$ for some $n$.  Rather than try to work with (or even properly define) nontrivial bundles, we will work with just a single patch $U=\mathbb R^n$.

Now, a common observation in mathematics is that things become easier when we work over the complex numbers. The set of complex numbers $\mathbb C =$\linebreak$\{a+bi|a,b\in\mathbb R\}$ is in correspondence with the set of points $\mathbb R^2=\{(a,b)|a,b\in\mathbb R\}$.\linebreak We will consider a base patch $U=\mathbb R^2$ and regard this as $U=\mathbb C$. Such patches cover a special class of base spaces called \emph{Riemann surfaces}, which are the typical setting for the theory of Higgs bundles.

\begin{figure}[ht]
\centering
\resizebox{0.75\textwidth}{!}{%
  \includegraphics{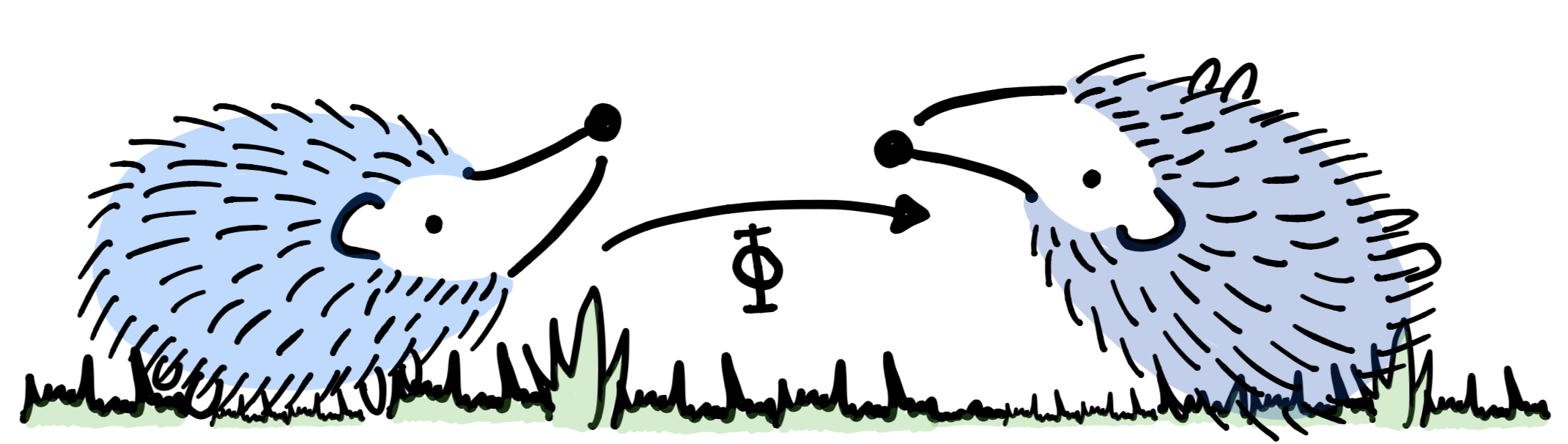}   
  }
\caption{A map  $\Phi$ which ``twists'' the hedgehog bundle.}\label{Lau2}
\end{figure}

A Higgs bundle is a vector bundle that also comes with a map $\Phi$, called the \emph{Higgs field}, which transforms the vector bundle in a certain way.  This is sometimes called ``twisting'' (see Figure \ref{Lau2}).

Because we have elected to work on a single patch of a Riemann surface, the map $\Phi$ can be represented in a particularly simple way.  If we have a complex rank-$r$ bundle on $U=\mathbb C$ (which means that a copy of $\mathbb C^r=\mathbb R^{2r}$ appears over every point of $U$), then $\Phi$ is an $r\times r$ matrix with polynomial entries.  In other words, the entries of $\Phi$ are polynomial functions of some degree $k$  taking the form $$a_kz^k+\cdots+a_1z+a_0,$$with $a_i\in\mathbb C$ for each $i$.  The parameter $z$ is a complex number that tells us where we are in the base $U$.  When we pick a particular value for $z$, the map $\Phi$ becomes an ordinary matrix of complex numbers and it acts on the copy of $\mathbb C^r$ over $z\in U$ by matrix multiplication.  Here, the degree of the polynomials can be regarded as the ``twisting'' effect.

Since $z$ encodes $U$ and since the size of $\Phi$ as a matrix reminds us that the bundle consists of copies of $\mathbb C^r$, we can go so far as to forget the patch and forget the $\mathbb C^r$'s altogether, and just focus on the map $\Phi$.  This is an extremely reductionist point of view, but it is useful for us that matrices with polynomial entries exhibit many of the interesting features of true Higgs bundles, since it uses only linear algebra   without many of the geometric complications.

\section{Moduli spaces of matrices with polynomial entries and spectral curves}

To start with the simplest example, we can take matrices with polynomial entries in which the entries are simply complex numbers.  These are, if one likes, the matrices with polynomial entries of degree zero.  We will work out the moduli space of complex $2\times2$ matrices.  We must first decide what it means for two matrices to be  equivalent.  The natural relation is that two matrices $A$ and $B$ are equivalent if they are \emph{similar}, meaning that there exists an invertible matrix $P$ in $\mbox{GL}(2,\mathbb{C})$ such that \[B=P^{-1}AP.\] 
For example consider the $2\times2$ matrices 
 \begin{eqnarray}A=\left(\begin{array}{cr}a&-b\\b&a\end{array}\right),  ~\hspace{0.8cm}~B=\left(\begin{array}{rc}a&b\\-b&a\end{array}\right),  ~\hspace{0.8cm}~ P=\left(\begin{array}{rc}-1&0\\0&1\end{array}\right).\nonumber\end{eqnarray}
 One can check that $B=P^{-1}AP$ and so $A$ and $B$ are equivalent. In fact, writing the complex numbers $a+bi$ as real $2\times2$ matrices $A$, the transformation $P$ corresponds to complex conjugation, $a+bi\mapsto a-bi$.

It turns out that any $2\times2$ matrix $A$ is equivalent to a matrix of one of the following three forms:
  \begin{eqnarray}D=\left(\begin{array}{cc}\lambda_1&0\\0&\lambda_2\end{array}\right),  ~\hspace{0.8cm}~D_0=\left(\begin{array}{cc}\lambda_1&0\\0&\lambda_1\end{array}\right),  ~\hspace{0.8cm}~ D_1=\left(\begin{array}{cc}\lambda_1&1\\0&\lambda_1\end{array}\right)\label{D}\end{eqnarray}\noindent where $\lambda_1\neq\lambda_2$.  The numbers $\lambda_1$ and $\lambda_2$ are simply the eigenvalues of $A$.  The fact that we can transform $A$ into an equivalent representative in which these numbers emerge as the only non-trivial data suggests that these values play a role in organizing the moduli space of the matrices.  To go back to the telephone book analogy, it appears that we can look up an equivalence class of similar $2\times2$ matrices by its eigenvalues.

When $\lambda_1\neq\lambda_2$, the pair $(\lambda_1,\lambda_2)$ corresponds to a single class of matrices, the one represented by $D$ in Equation~\eqref{D}.  Something a little different happens when we encounter a pair of the form $(\lambda_1,\lambda_1)$.  Here, there are exactly two classes that share this ``telephone number'': the class of matrices similar to $D_0$ and the class of matrices similar to $D_1$ in Equation~\eqref{D}.  If it bothers you that most eigenvalue pairs correspond to only one class of matrices but some pairs correspond to two classes, then you are not alone.  This is a fundamental problem in moduli theory and algebraic geometers solve this by throwing away the extra bothersome points.  In this case, we will choose to discard the classes represented by $D_1$ in Equation~\eqref{D}.  This act of ``throwing away'' is what algebraic geometers refer to as a \emph{imposing a stability condition}.  In linear algebra terms, we are throwing away the non-diagonalizable matrices from the set of all $2\times2$ matrices.  

From this, it is tempting to conclude that the space $\mathbb C^2$, this is, the space of pairs $(\lambda_1,\lambda_2)$, is the moduli space of diagonalizable $2\times2$ matrices under similarity.  Indeed, for each point $(\lambda_1,\lambda_2)\in\mathbb C^2$, regardless of whether $\lambda_1=\lambda_2$ or not, there is a single class of such matrices.  At this point, one would do well to recognize the following subtlety: the matrices$$\left(\begin{array}{cc}\lambda_1&0\\0&\lambda_2\end{array}\right)~~\mbox{ and }~~\left(\begin{array}{cc}\lambda_2&0\\0&\lambda_1\end{array}\right)$$are equivalent under similarity.  In other words, the pairs $(\lambda_1,\lambda_2)$ and $(\lambda_2,\lambda_1)$ represent the same class. The solution is to make no distinction between any two pairs that differ only by order. In group-theoretic terms, this is the same as saying that the moduli space is the quotient of $\mathbb C^2=\mathbb C\times\mathbb C$ by $S_2$, the symmetric group on two letters.  The beautiful thing is that the space that results from collecting pairs in this way is again in bijection with $\mathbb C^2$ as a set,\footnote{Further still, they are \emph{homeomorphic} as topological spaces. What is lost after quotienting is the vector space structure on $\mathbb C^2$, and so algebraic geometers would prefer to refer to the moduli space as the \emph{affine space} $\mathbb A^2$.} and so we may indeed conclude that the moduli space is $\mathbb C^2$. Viewed from the theory of Higgs bundles, the space $\mathbb C^2$ plays the role of what we call the \emph{Hitchin base}, while the single class of matrices corresponding to each point in $\mathbb C^2$ is the \emph{Hitchin fibre}.  

\begin{figure}[ht]
\centering
\resizebox{0.35\textwidth}{!}{%
  \includegraphics{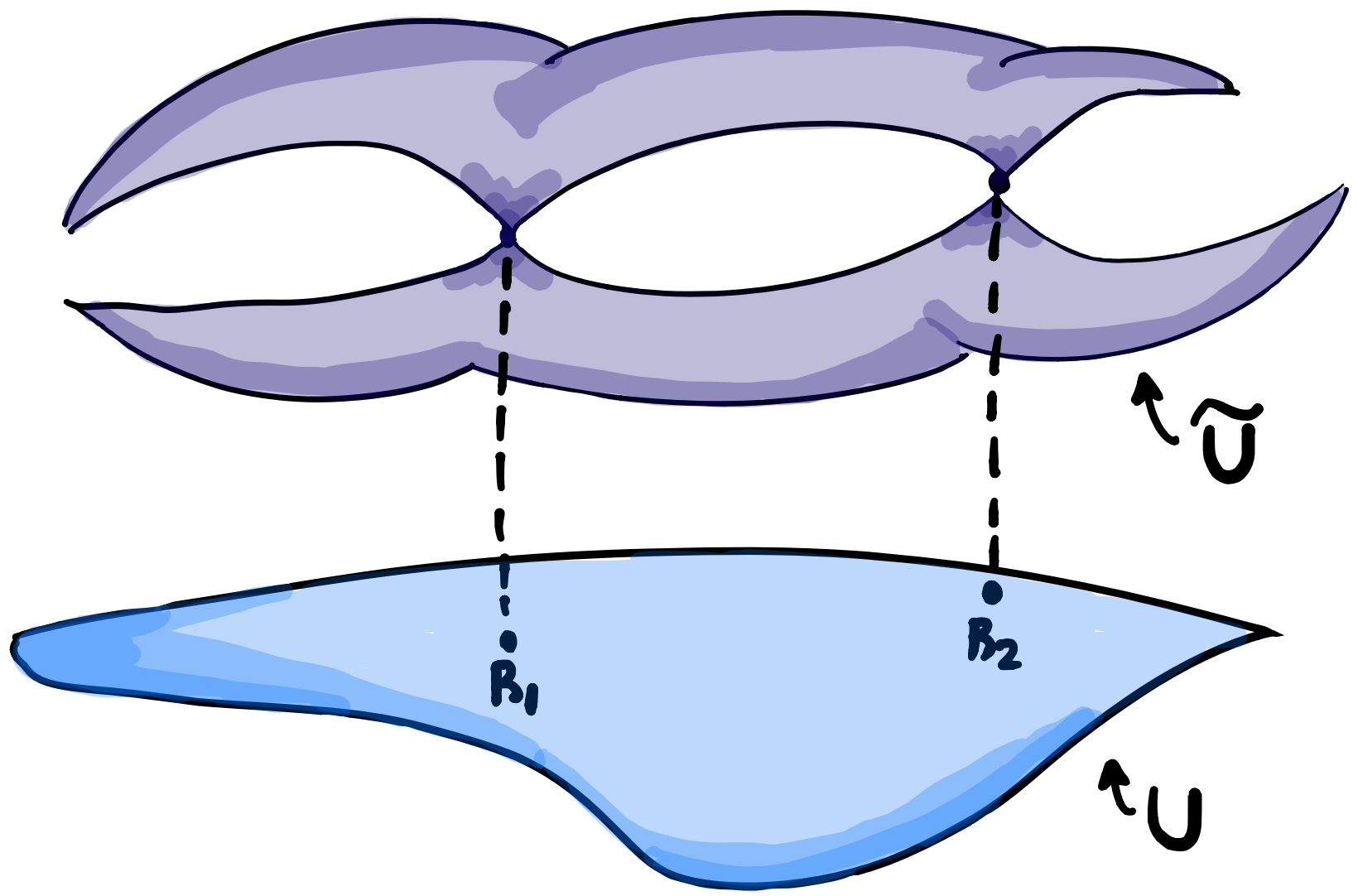}   
  }
\caption{A spectral curve with $R_1$ and $R_2$ corresponding to two branch points.}\label{Lau4}
\end{figure}

To go further, it is important to remember that the two eigenvalues $\lambda_1$ and $\lambda_2$ in the preceding example emerge as solutions to a degree $2$ polynomial equation,$$\lambda^2-(\mbox{tr}\,A)\lambda+(\det A)=0,$$the \emph{characteristic equation} of $A$, where $\mbox{tr}$ and $\det$ are the trace and determinant, respectively.\footnote{For each characteristic equation, there is exactly one pair of (unordered) eigenvalues that solve the equation, and therefore exactly one diagonalizable class of $2\times2$ matrices.  This is another way to see that the moduli space of $2\times2$ matrices is $\mathbb C^2$, since it takes two numbers, $\mbox{tr}(A)$ and $\det(A)$, to write down each characteristic equation.  In this model, the order of the two numbers does make a difference --- swapping a trace with a determinant will  change the equation and hence the solutions!}  When $\Phi$ is a $2\times2$ polynomial-valued matrix of degree larger than zero, both $\mbox{tr}\,\Phi$ and $\det \Phi$ will themselves be polynomials in $z$.  This means that for each point $z$ in the patch $U$, there will be a pair of eigenvalues $(\lambda_1(z),\lambda_2(z))$.  This forms a new patch $\widetilde U$ that, for almost every $z\in U$, contains two points that project to $z$.  This is what geometers call a \emph{branched double cover}.  The word ``branched'' accounts for the possibility that, for some $z$, there will be equal eigenvalues where the two sheets of $\widetilde U$ will come together as in Figure \ref{Lau4}.
This new patch $\widetilde U$ is a local model for what geometers call a \emph{spectral curve} --- a space whose points are the eigenvalues of some matrix or family of matrices.  

\begin{figure}[ht]
\centering
\resizebox{0.55\textwidth}{!}{%
  \includegraphics{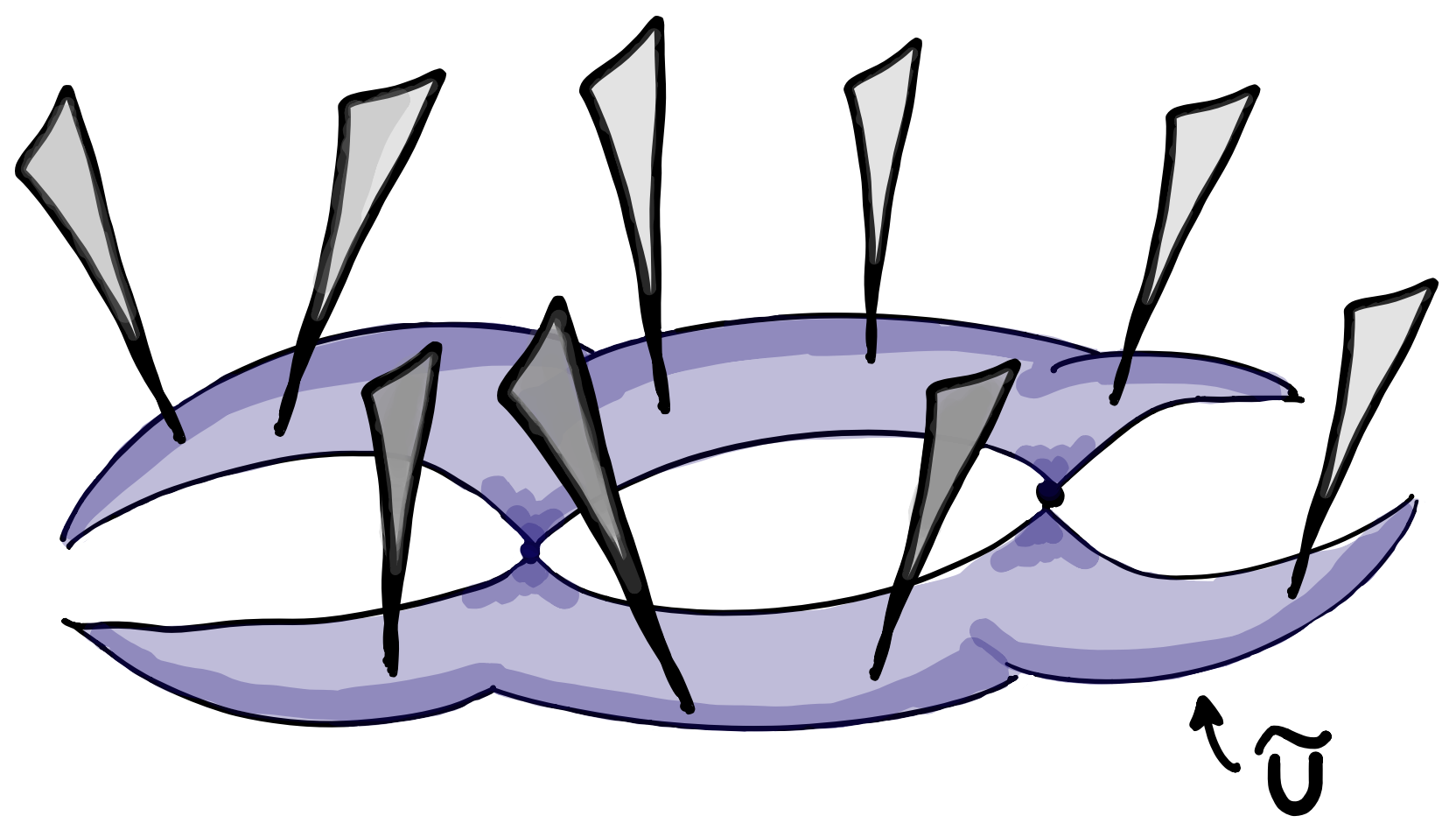}   
  }
\caption{The data defining the vector bundle of a Higgs bundle, seen as a line bundle on the double cover of the base space. }\label{Lau5}
\end{figure}

If we populate $\widetilde U$ with the structure of a complex line bundle, as illustrated in Figure \ref{Lau5}, then for most points $z\in U$ there will be two complex lines situated directly above them in $\widetilde U$.  We can physically transport these lines down to $z$ --- this operation is called the \emph{pushforward} --- where they will span a copy of $\mathbb C^2$ at $z$.  Furthermore, we can read off the two eigenvalues $\lambda_1(z)$ and $\lambda_2(z)$ above $z$ and construct a diagonal matrix $\Phi_z$ from them.  It is possible to formulate $\Phi_z$ in a consistent way even at points $z$ where $\widetilde U$ is branched.  Having done this for all $z\in U$, we will have produced a rank-$2$ complex vector bundle and a corresponding polynomial-valued matrix $\Phi$ on $U$.  It is possible to go back and forth between these types of objects: matrices with polynomial entries on $U$ and complex line bundles on $\widetilde U$.  This is called the \emph{spectral correspondence} and it extends even to non-trivial Higgs bundles on compact Riemann surfaces.  A complete account of this appears in the work of Hitchin \cite{Hi2}.

\section{A hint of the full picture}

\begin{figure}[ht]
\centering
\resizebox{0.8\textwidth}{!}{%
  \includegraphics{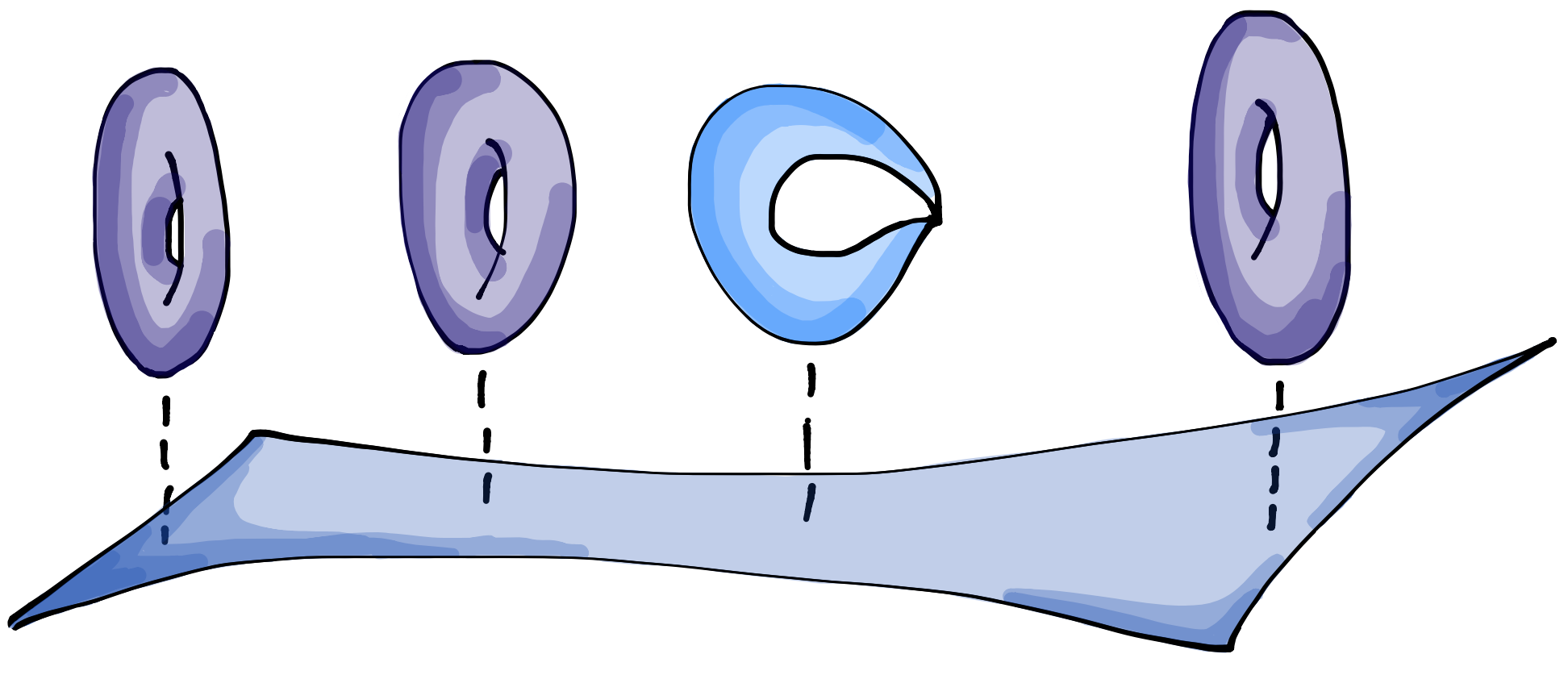}   
  }
\caption{A torus fibration with singular fibres.}\label{Lau6}
\end{figure}

In the full geometric picture on a Riemann surface $X$, each Higgs bundle on $X$\linebreak will determine a spectral curve $\widetilde X$ that forms a branched cover of $X$.  The spectral correspondence replaces Higgs bundles on $X$ with line bundles on the corresponding $\widetilde X$ and the correspondence respects equivalence, which is a generalization of similarity of matrices.  If two Higgs bundles are equivalent, then they determine the same $\widetilde X$ and their line bundles on $\widetilde X$ are also equivalent.  The set of equivalence classes of Higgs bundles with the same spectral curve $\widetilde X$\linebreak is in correspondence with a set of line bundles\footnote{We should be discussing \emph{holomorphic Higgs bundles} and \emph{holomorphic line bundles} at this point, but we wish to maintain the decidedly very informal tone of this note for the lay reader.} on $\widetilde X$, and these form a geometric torus (of some dimension).  The collection of all spectral curves is the Hitchin base while the respective tori are the Hitchin fibres, see Figure \ref{Lau6}.\linebreak As with ordinary matrices, decisions have to be made regarding throwing away ``unstable'' Higgs bundles, in order to build a   moduli space.  In the case of ordinary matrices, the spectral curves become ordinary eigenvalues and line bundles become ordinary eigenspaces, of which there is no choice up to equivalence,  recovering the picture discussed earlier in which the Hitchin base consists of eigenvalues and the Hitchin fibres are just single points.
\enlargethispage{\baselineskip}

\pagebreak

\emph{Torus fibrations}, which are geometric structures consisting of a torus over each point in some other space, appear frequently in mathematics and physics: in the theory of integrable systems, in mirror symmetry, and in representation theory.  This universality evokes natural questions about how the moduli space of Higgs bundles might connect with these areas.  It happens to be the case that Higgs bundles play a crucial role in problems central to all three subjects.

If we have kept the attention of the reader, we invite them to spend some time with introductory surveys on Higgs bundles such as \cite{Hi3, steve, lau} and, from there, examine some more specialized literature on the rich interplay between Higgs bundles and various areas of mathematics and physics.


\begin{imagecredits}
\item[] All figures were made by the authors. 
\end{imagecredits}

\pagebreak


\end{document}